\documentclass[12pt]{article}
\usepackage{graphicx}
\usepackage{amsmath,amsthm,amssymb,enumerate}%, esint}
\usepackage{euscript,mathrsfs}
\usepackage[left=2cm,right=2cm,top=3.5cm,bottom=3.5cm]{geometry}
\usepackage{color}
\allowdisplaybreaks

\usepackage{soul}

\catcode`\@=11 \@addtoreset{equation}{section}

\catcode`\@=12

\allowdisplaybreaks

\newtheorem{Theorem}{Theorem}[section]
\newtheorem{Proposition}[Theorem]{Proposition}
\newtheorem{Lemma}[Theorem]{Lemma}
\newtheorem{Corollary}[Theorem]{Corollary}

\theoremstyle{definition}
\newtheorem{Definition}[Theorem]{Definition}

\newtheorem{Remark}[Theorem]{Remark}

\newcommand{\bTheorem}[1]{
\begin{Theorem} \label{T#1} }
\newcommand{\eT}{\end{Theorem}}

\newcommand{\bProposition}[1]{
\begin{Proposition} \label{P#1}}
\newcommand{\eP}{\end{Proposition}}

\newcommand{\bLemma}[1]{
\begin{Lemma} \label{L#1} }
\newcommand{\eL}{\end{Lemma}}

\newcommand{\bCorollary}[1]{
\begin{Corollary} \label{C#1} }
\newcommand{\eC}{\end{Corollary}}

\newcommand{\bRemark}[1]{
\begin{Remark} \label{R#1} }
\newcommand{\eR}{\end{Remark}}

\newcommand{\bDefinition}[1]{
\begin{Definition} \label{D#1} }
\newcommand{\eD}{\end{Definition}}

\newcommand{\br}{\nonumber \\}

\newcommand{\Del}{\Delta_x}

\newcommand{\vrB}{\vr_B}

\newcommand{\Ds}{\mathbb{D}_x}

\newcommand{\vuB}{\vc{u}_B}

\newcommand{\tvm}{\widetilde{\vc{m}}}
\newcommand{\tS}{\widetilde{S}}
\newcommand{\bfphi}{\boldsymbol{\varphi}}

\newcommand{\bFormula}[1]{
\begin{equation} \label{#1}}
\newcommand{\eF}{\end{equation}}

\newcommand{\Ov}[1]{\overline{#1}}

\newcommand{\DC}{C^\infty_c}

\newcommand{\aleq}{\stackrel{<}{\sim}}

\newcommand{\ageq}{\stackrel{>}{\sim}}

\newcommand{\vr}{\varrho}

\newcommand{\tvr}{\tilde \vr}
\newcommand{\tvu}{{\tilde \vu}}
\newcommand{\tvt}{\tilde \vt}

\newcommand{\vt}{\vartheta}
\newcommand{\vu}{\vc{u}}
\newcommand{\vm}{\vc{m}}

\newcommand{\vtB}{\vartheta_B}
\newcommand{\vc}[1]{{\bf #1}}

\newcommand{\Div}{{\rm div}_x}
\newcommand{\Grad}{\nabla_x}

\newcommand{\dx}{\,{\rm d} {x}}

\newcommand{\dt}{\,{\rm d} t }

\newcommand{\intO}[1]{\int_{\Omega} #1 \ \dx}

\newcommand{\D}{{\rm d}}

\newcommand{\ep}{\varepsilon}

\newcommand{\intgi}[1]{\int_{\partial \Omega} #1 \ [ \vuB \cdot \vc{n} ]^- \D \sigma_x}
\newcommand{\intgo}[1]{\int_{\partial \Omega} #1 \ [ \vuB \cdot \vc{n} ]^+ \D \sigma_x}

\def\softd{{\leavevmode\setbox1=\hbox{d}%
          \hbox to 1.05\wd1{d\kern-0.4ex{\char039}\hss}}}
\definecolor{Cgrey}{rgb}{0.85,0.85,0.85}
\definecolor{Cblue}{rgb}{0.50,0.85,0.85}
\definecolor{Cred}{rgb}{1,0,0}
\definecolor{fancy}{rgb}{0.10,0.85,0.10}

\newcommand\Cbox[2]{%
    \newbox\contentbox%
    \newbox\bkgdbox%
    \setbox\contentbox\hbox to \hsize{%
        \vtop{
            \kern\columnsep
            \hbox to \hsize{%
                \kern\columnsep%
                \advance\hsize by -2\columnsep%
                \setlength{\textwidth}{\hsize}%
                \vbox{
                    \parskip=\baselineskip
                    \parindent=0bp
                    #2
                }%
                \kern\columnsep%
            }%
            \kern\columnsep%
        }%
    }%
    \setbox\bkgdbox\vbox{
        \color{#1}
        \hrule width  \wd\contentbox %
               height \ht\contentbox %
               depth  \dp\contentbox
        \color{black}
    }%
    \wd\bkgdbox=0bp%
    \vbox{\hbox to \hsize{\box\bkgdbox\box\contentbox}}%
    \vskip\baselineskip%
}

\date{}

\begin{document}

%%%%%%%%%%%%%%%%%%%%%%%%%%%%%%%%

\title{Navier--Stokes--Fourier system with Dirichlet boundary conditions}

\author{ Nilasis Chaudhuri \footnotemark[1] \and Eduard Feireisl
\thanks{The work of N.C. and E.F. was supported by the
Czech Sciences Foundation (GA\v CR), Grant Agreement
18--05974S. The Institute of Mathematics of the Academy of Sciences of
the Czech Republic is supported by RVO:67985840.}
}

%\date{\today}

\maketitle

\centerline{Institute f{\"u}r Mathematik, Technische Universit{\"a}t, Berlin}

\centerline{Stra\ss e des 17. Juni 136, D -- 10623 Berlin, Germany.}
\centerline{chaudhuri@math.tu-berlin.de}
\centerline {and}
\centerline{Institute of Mathematics of the Academy of Sciences of the Czech Republic;}
\centerline{\v Zitn\' a 25, CZ-115 67 Praha 1, Czech Republic}

\centerline{feireisl@math.cas.cz}

\begin{abstract}
	
	We consider the Navier--Stokes--Fourier system describing the motion of a compressible, viscous, and heat conducting fluid in a bounded domain 
	$\Omega \subset R^d$, $d =2,3$, with general non-homogeneous Dirichlet boundary conditions for the velocity and the absolute temperature, with the associated 
	boundary conditions for the density on the inflow part. We introduce a new concept of weak solution based on the satisfaction of the entropy inequality together with a 
	balance law for the ballistic energy. We show the weak--strong uniqueness principle as well as the existence of global--in--time solutions.

\end{abstract}


%{\bf MSC:}
\bigskip

\centerline{\it To the memory of Andro Mikeli\' c}

%\tableofcontents

\section{Introduction}
\label{I}

The time evolution of the \emph{density} $\vr = \vr(t,x)$, the bulk \emph{velocity} $\vu = \vu(t,x)$ and the (absolute) \emph{temperature} $\vt= \vt(t,x)$ of a general compressible, viscous, and heat conducting fluid 
is governed by the following system of field equations:
\begin{align} 
	\partial_t \vr + \Div (\vr \vu) &= 0, \br
	\partial_t (\vr \vu) + \Div (\vr \vu \otimes \vu) + \Grad p(\vr, \vt) &= \Div \mathbb{S}+ \vr \vc{g}, \ \br
	\partial_t (\vr e(\vr, \vt)) + \Div (\vr e(\vr, \vt) \vu ) + \Div \vc{q} &= \mathbb{S} : \Ds \vu - p(\vr, \vt) \Div \vu,\ \Ds \vu \equiv \frac{1}{2} (\Grad \vu + \Grad^t \vu).
	\label{i1}
	\end{align}
Here $p = p(\vr, \vt)$ is the pressure related to the internal energy $e(\vr, \vt)$ through \emph{Gibbs' equation} 
\begin{equation} \label{i2}
	\vt D s = D e + p D \left( \frac{1}{\vr} \right), 
	\end{equation}
where the quantity $s = s(\vr, \vt)$ is called \emph{entropy}. We suppose that the fluid is \emph{Newtonian}, meaning the \emph{viscous stress tensor} 
$\mathbb{S}$ satisfies 
\begin{equation} \label{i3}
	\mathbb{S} (\vt, \Ds \vu) = \mu(\vt) \left( \Grad \vu + \Grad^t \vu - \frac{2}{d} \Div \vu \mathbb{I} \right) + \eta(\vt) \Div \vu \mathbb{I},
\end{equation}
with the \emph{shear viscosity} coefficient $\mu$ and the \emph{bulk viscosity} coefficient $\eta$. Similarly, we suppose that the \emph{heat flux} $\vc{q}$ is determined by 
\emph{Fourier's law}, 
\begin{equation} \label{i4}
	\vc{q} (\vt, \Grad \vt) = - \kappa (\vt) \Grad \vt,
	\end{equation}
with the \emph{heat conductivity} coefficient $\kappa$. The symbol $\vc{g} = \vc{g}(t,x)$ stands for a given external force acting on the fluid. The system of equations \eqref{i1}, together 
with the constitutive relations \eqref{i3}, \eqref{i4} is called \emph{Navier--Stokes--Fourier} system.
 
We suppose the fluid is confined to a bounded regular domain $\Omega \subset R^d$, $d = 2,3$. Our main goal is to discuss solvability of the initial--boundary value problem for the Navier--Stokes--Fourier system 
endowed with the \emph{non--homogeneous Dirichlet boundary conditions} 
\begin{align}
	\vu|_{\partial \Omega} &= \vuB , \br
	\vt|_{\partial \Omega} &= \vtB 
	\label{i5}
	\end{align}
Accordingly, the \emph{inflow boundary conditions} for the density must be prescribed, 
\begin{equation} \label{i6} 
	\vr|_{\Gamma_{\rm in}} = \vrB,\ \Gamma_{\rm in} \equiv \left\{ (t,x) \ \Big| t \in [0,T],\ x \in \partial \Omega,\ \vuB(t,x) \cdot \vc{n}(x) < 0 \right\},
	\end{equation}
where $\vc{n}$ denotes the outer normal vector to $\partial \Omega$.

The problem \eqref{i1}--\eqref{i6} is well posed \emph{locally in time} for smooth initial data in the class of smooth solutions, see e.g. Valli and Zajaczkowski \cite{VAZA}. Our objective is global--in--time solvability in the class of weak solutions. The overwhelming majority of the available literature on the weak solvability of the Navier--Stokes--Fourier system is devoted to conservative or space periodic boundary conditions, 
where the total energy of the system is conserved, see Lions \cite[Chapter 8]{LI4}, Bresch and Desjardins \cite{BRDE}, \cite{BRDE1}, Bresch and Jabin \cite{BreJab} or the monographs \cite{EF70}, 
\cite{FeNo6A}. Recently, the approach of \cite{FeNo6A} has been extended to the open system with general in/out flow boundary conditions, see \cite{FeiNov20}. Still, the choice of boundary conditions 
admissible in \cite{FeiNov20} requires the control of the internal (heat) energy flux $\vc{q}$ on $\partial \Omega$. Thus the important class of problems, where the boundary temperature $\vtB$ is given, among 
which the well known Rayleigh-B\' enard problem, is not covered by the theory developed in \cite{FeiNov20}. The goal of the present paper is to fill this gap.

The principal and well known difficulty of nonlinear systems of fluid mechanics is that the {\it a priori} bounds based on the energy estimates are not strong enough to render certain quantities, specifically the source term 
$\mathbb{S} : \Ds \vu - p \Div \vu$, equi--integrable. The approach used in \cite{FeNo6A} replaces the internal energy equation by the entropy \emph{inequality} 
\begin{equation} \label{i7}
\partial_t (\vr s(\vr, \vt)) + \Div (\vr s(\vr \vt) \vu) + \Div \left( \frac{ \vc{q}(\vt, \Grad \vt ) }{\vt} \right) \geq \frac{1}{\vt} 
\left( \mathbb{S} (\vt, \Ds \vu ) : \Ds \vu - \frac{\vc{q}(\vt, \Grad \vt) \cdot \Grad \vt }{\vt} \right)
\end{equation}
supplemented by the total energy balance 
\begin{equation} \label{i8}
\frac{\D }{\dt} \intO{ \left[ \frac{1}{2} \vr |\vu|^2 + \vr e(\vr, \vt) \right] } = \intO{ \vr \vc{g} \cdot \vu }.
\end{equation}
Of course, the energy balance \eqref{i8} holds only for energetically insulated systems, where $\vuB = 0$, $\vc{q} \cdot \vc{n}|_{\partial \Omega} = 0$. The main idea to accommodate the Dirichlet boundary conditions for the temperature is to replace \eqref{i8} by a similar relation for the \emph{ballistic energy} 
\[
\intO{ \left[ \frac{1}{2} \vr |\vu - \vuB|^2 + \vr e(\vr, \vt) - \tvt \vr s(\vr, \vt) \right] },  
\]
where $\tvt$ is an arbitrary smooth function satisfying the boundary conditions $\tvt|_{\partial \Omega} = \vtB$. Here the term  \emph{ballistic} is motivated by Ericksen \cite{Eri}, where the quantity 
\[
\vr e - \vtB \vr s 
\]
is called \emph{ballistic free energy}. The advantage of working with the ballistic energy rather than the total energy is that the time evolution of the former can be described only 
in terms of the boundary values $\vuB$, $\vtB$, $\vrB$. 

The weak formulation based on the balance of entropy and ballistic free energy gives rise to the concept of \emph{weak solution} enjoying the following properties:
\begin{itemize}
	\item {\bf Global existence.} The weak solutions exist globally in time for any physically admissible choice of the initial/boundary data.
	\item {\bf Weak-strong uniqueness.} A weak solution and the strong solution corresponding to the same initial/boundary data coincide as long as the strong solution exists.
	\item {\bf Compatibility.} Any regular weak solution is a strong solution.
	
	\end{itemize}

The paper is organized as follows. In Section \ref{w}, we collect the available preliminary material and introduce the concept of weak solution. In Section \ref{WS}, we prove the weak--strong uniqueness 
property. Finally, in Section \ref{E}, we show the existence of weak solutions on an arbitrary time interval.

\section{Weak formulation}
\label{w}

Our goal is to derive a suitable \emph{weak formulation} for the Navier--Stokes--Fourier system supplemented with the boundary conditions \eqref{i5}, \eqref{i6}.

\subsection{Equation of continuity}

We say that $\vr$, $\vu$ satisfy the equation of continuity (\ref{i1}$)_1$ if 
the integral identity
\begin{align}
	\int_0^\tau &\intO{ \Big[ \vr \partial_t \varphi + \vr \vu \cdot \Grad \varphi \Big]} \dt \br 
	&= \intgi{ \varphi \vrB } + \intgo{ \varphi \vr} + \left[ \intO{ \varphi \vr } \right]_{t = 0}^{t = \tau}
	\label{w1}
\end{align}
holds for any $0 \leq \tau \leq T$, $\varphi \in C^1([0,T] \times \Ov{\Omega})$. Here and hereafter, we use the notation 
\[
[v]^+ = \max \{ v,0 \},\ [v]^- = \min \{ v , 0 \}.
\]
The quantity on the right--hand side should be interpreted as the normal trace of the 
$(d+1)-$dimensional vector field $[\vr, \vu]$ on the space--time cylinder $(0, \tau) \times \Omega$, cf. Chen, Torres, Ziemer \cite{ChToZi}. We tacitly assume that 
all integrals exist finite.

\subsection{Momentum equation}

We suppose that the velocity admits a Sobolev trace $\vuB$ on $\partial \Omega$, specifically, 
\begin{equation} \label{w2}
	(\vu - \vuB) \in L^r(0,T; W^{1,r}_0(\Omega; R^d)) \ \mbox{for some}\ r > 1,
\end{equation} 
where $\vuB = \vuB(t,x)$ has been extended to all $(t,x) \in R^{d+1}$. The momentum equation (\ref{i1}$)_2$ is then interpreted in the following sense:
\begin{align}
	\int_0^\tau &\intO{ \Big[ \vr \vu \cdot \partial_t \bfphi + \vr \vu \otimes \vu : \Grad \bfphi + p(\vr, \vt) \Div \bfphi \Big] } \dt \br &= 
	\int_0^\tau \intO{ \Big[ \mathbb{S} (\vt, \Ds \vu) : \Ds \bfphi - \vr \vc{g} \cdot \bfphi \Big] } \dt + \left[ \intO{ \vr \vu \cdot \bfphi } \right]_{t=0}^{t = \tau}  
	\label{w3}
\end{align}
for any $0 \leq \tau \leq T$, $\bfphi \in C^1_c([0,T] \times \Omega; R^d)$.

\subsection{Entropy inequality}

As pointed out in the introduction, our strategy is to ``replace'' the internal energy balance (\ref{i1}$)_3$ by the entropy inequality \eqref{i7}. To comply with the Dirichlet boundary condition 
for the temperature, we suppose 
\begin{equation} \label{w4}
	(\vt - \vtB) \in L^2(0,T; W^{1,2}_0(\Omega)),
	\end{equation}
where $\vtB = \vtB(t,x)$ has been extended as a smooth strictly positive function on $R^{d+1}$. The weak formulation of \eqref{i7} reads
\begin{align}
	&- \int_0^T \intO{ 
		\left[ \vr s \partial_t \varphi + \vr s \vu \cdot \Grad \varphi + \frac{\vc{q}}{\vt} \cdot \Grad \varphi \right] } \dt \br &\geq 
	\int_0^T \intO{ \frac{\varphi}{\vt} \left( \mathbb{S}(\vt, \Ds \vu ): \Ds \vu - \frac{\vc{q}(\vt, \Grad \vt) \cdot \Grad \vt }{\vt} \right) } \dt 
	+ \intO{ (\vr s)(0, \cdot) \varphi }
	\label{w5}
\end{align}
for any $\varphi \in C^1_c([0,T) \times \Omega)$, $\varphi \geq 0$. Unlike the density $\vr$ or the momentum $\vr \vu$ the total entropy $\vr s$ is in general not weakly continuous 
in the time variable. Nevertheless, it follows from \eqref{w5} that 
\[
t \in [0,T] \mapsto \intO{ (\vr s) (t, \cdot) \phi },\ \phi \in C^1_c(\Omega),\ \phi \geq 0,
\]
can be written as a sum of a non--decreasing and a continuous function. Consequently, setting 
\[
\intO{ \vr s (\tau, \cdot) \phi } = \lim_{\delta \searrow 0} \frac{1}{\delta} \int_{\tau - \delta}^\tau \intO{ \vr s \phi },\ 0 < \tau \leq T,
\]
we may consider the total entropy $\vr s$ is weakly c\` agl\` ad (left continuous with the limit from the right) function of $t \in [0,T]$ with the convention that its value at $t=0$ is given by the initial data.
Accordingly, the weak formulation \eqref{w5} can be written in the form
\begin{align}
	&\left[ \intO{ \vr s \varphi } \right]_{t = \tau_1}^{t = \tau_2} - \int_{\tau_1}^{\tau_2} \intO{ 
		\left[ \vr s \partial_t \varphi + \vr s \vu \cdot \Grad \varphi + \frac{\vc{q}}{\vt} \cdot \Grad \varphi \right] } \dt \br &\geq 
	\int_{\tau_1}^{\tau_2} \intO{ \frac{\varphi}{\vt} \left( \mathbb{S}(\vt, \Ds \vu ): \Ds \vu - \frac{\vc{q}(\vt, \Grad \vt) \cdot \Grad \vt }{\vt} \right) } \dt 
	\label{w6}
\end{align}
for any $0 \leq \tau_1 < \tau_2 \leq T$, and any $\varphi \in C^1_c([0,T] \times \Omega)$, $\varphi \geq 0$.

\subsection{Ballistic energy}

Assuming for a moment that all quantities in question are smooth,  we
multiply the momentum equation on $(\vu - \vuB)$ and integrate by parts over $\Omega$. Similarly, 
we integrate the internal energy equation (\ref{i1}$)_3$. Summing up the results, we obtain the total energy balance
\begin{align}  
	\frac{\D }{\dt} &\intO{ \left[ \frac{1}{2} \vr |\vu - \vuB|^2 + \vr e \right] } + 
	\intgi{ \vrB e(\vrB, \vtB) } + \intgo{ \vr e(\vr, \vtB) } \br &+ 
	\int_{\partial \Omega} \vc{q} \cdot \vc{n} \ \D \sigma_x  
	= - \intO{ \Big[ \vr \vu \otimes \vu + p \mathbb{I} - \mathbb{S} \Big] : \Ds \vuB } + \frac{1}{2} \intO{ \vr \vu \cdot \Grad |\vuB|^2 } \br &\ \ \ + 
	\intO{ \vr (\vu - \vuB)\cdot (\vc{g} - \partial_t \vuB) } 
\label{w7}
	\end{align}
Unfortunately, relation \eqref{w7} cannot be used in the weak formulation as 
we do not control the boundary integral 
\[
\int_{\partial \Omega} \vc{q} \cdot \vc{n} \ \D \sigma_x .
\]

Consider $\tvt \in C^1([0,T] \times \Ov{\Omega})$, $\tvt > 0$, $\tvt|_{\partial \Omega} = \vtB$. Multiplying the entropy inequality \eqref{i7} by $\tvt$ and integrating by parts we obtain 
\begin{align} 
	- \frac{\D }{\dt} &\intO{ \tvt \vr s } - \int_{\partial \Omega} \vr \vtB s \vuB \cdot \vc{n} \ \D \sigma_x  - \int_{\partial \Omega} \vc{q} \cdot \vc{n} \D \sigma_x \leq 
	- \intO{ \frac{\tvt}{\vt}	 \left( \mathbb{S} : \Ds \vu - \frac{\vc{q} \cdot \Grad \vt }{\vt} \right) } \br 
	&- \intO{ \left[ \vr s \left( \partial_t \tvt + \vu \cdot \Grad \tvt \right) + \frac{\vc{q}}{\vt} \cdot \Grad \tvt \right] }
	\label{w8}
\end{align}

Summing up \eqref{w7}, \eqref{w8} we get the \emph{ballistic energy inequality}
\begin{align}  
	\frac{\D }{\dt} &\intO{ \left[ \frac{1}{2} \vr |\vu - \vuB|^2 + \vr e - \tvt \vr s \right] } \br &+ 
	\intgi{ \Big[ \vrB e(\vrB, \vtB) - \vtB \vrB s(\vrB, \vtB) \Big] }  \br &+ 
	\intgo{ \Big[ \vr e(\vr, \vtB) - \vtB \vr s(\vr, \vtB) \Big] } \br &+ 
	\intO{ \frac{\tvt}{\vt}	 \left( \mathbb{S} : \Ds \vu - \frac{\vc{q} \cdot \Grad \vt }{\vt} \right) } \br
	&\leq  - \intO{ \Big[ \vr \vu \otimes \vu + p \mathbb{I} - \mathbb{S} \Big] : \Ds \vuB } + \frac{1}{2} \intO{ \vr \vu \cdot \Grad |\vuB|^2 } \br &\ \ \ + 
	\intO{ \vr (\vu - \vuB)\cdot (\vc{g} - \partial_t \vuB) } \br 
	&\ \ \ - \intO{ \left[ \vr s \left( \partial_t \tvt + \vu \cdot \Grad \tvt \right) + \frac{\vc{q}}{\vt} \cdot \Grad \tvt \right] }.
	\label{w9}
\end{align}
Note carefully that \eqref{w9} does not contain the boundary heat flux and therefore it is suitable to be included in the weak formulation of the problem. Accordingly, 
we require 
\begin{align}  
	- \psi(0) & \intO{ \left( \frac{1}{2} \vr |\vu - \vuB|^2 + \vr e - \tvt \vr s \right) }  \br &- 
	\int_0^T \partial_t \psi \intO{ \left( \frac{1}{2} \vr |\vu - \vuB|^2 + \vr e - \tvt \vr s \right) } \dt  \br &+ 
	\int_0^T \psi \intgi{ \Big[ \vrB e(\vrB, \vtB) - \vtB \vrB s(\vrB, \vtB) \Big] } \dt \br &+ 
	\int_0^T \psi \intgo{ \Big[ \vr e(\vr, \vtB) - \vtB \vr s(\vr, \vtB) \Big] } \dt \br &+ 
	\int_0^T \psi \intO{ \frac{\tvt}{\vt}	 \left( \mathbb{S} : \Ds \vu - \frac{\vc{q} \cdot \Grad \vt }{\vt} \right) } \dt \br
	&\leq  - \int_0^T \psi \intO{ \Big[ \vr \vu \otimes \vu + p \mathbb{I} - \mathbb{S} \Big] : \Ds \vuB } \dt + \int_0^T \psi \frac{1}{2} \intO{ \vr \vu \cdot \Grad |\vuB|^2 } \dt \br &\ \ \ + 
	\int_0^T \psi \intO{ \vr (\vu - \vuB)\cdot (\vc{g} - \partial_t \vuB) } \dt \br 
	&\ \ \ - \int_0^T \psi \intO{ \left[ \vr s \left( \partial_t \tvt + \vu \cdot \Grad \tvt \right) + \frac{\vc{q}}{\vt} \cdot \Grad \tvt \right] } \dt
	\label{w10}
\end{align}
for any $\tvt \in C^1([0,T] \times \Ov{\Omega})$, $\tvt > 0$, $\tvt|_{\partial \Omega} = \vtB$, and any $\psi \in C^1_c [0,T)$, $\psi \geq 0$. By a proper choice of the test function $\psi$, it is a routine matter to deduce the integrated from of \eqref{w10}, namely 
\begin{align}  
	 & \left[ \intO{ \left( \frac{1}{2} \vr |\vu - \vuB|^2 + \vr e - \tvt \vr s \right) } \right]_{t = 0}^{t = \tau} \br & 
	\int_0^\tau \intgi{ \Big[ \vrB e(\vrB, \vtB) - \vtB \vrB s(\vrB, \vtB) \Big] } \dt \br &+ 
	\int_0^\tau \intgo{ \Big[ \vr e(\vr, \vtB) - \vtB \vr s(\vr, \vtB) \Big] } \dt \br &+ 
	\int_0^\tau \intO{ \frac{\tvt}{\vt}	 \left( \mathbb{S} : \Ds \vu - \frac{\vc{q} \cdot \Grad \vt }{\vt} \right) } \dt \br
	&\leq  - \int_0^\tau  \intO{ \Big[ \vr \vu \otimes \vu + p \mathbb{I} - \mathbb{S} \Big] : \Ds \vuB } \dt + \frac{1}{2} \int_0^\tau \intO{ \vr \vu \cdot \Grad |\vuB|^2 } \dt \br &\ \ \ + 
	\int_0^\tau \intO{ \vr (\vu - \vuB)\cdot (\vc{g} - \partial_t \vuB) } \dt \br 
	&\ \ \ - \int_0^\tau \intO{ \left[ \vr s \left( \partial_t \tvt + \vu \cdot \Grad \tvt \right) + \frac{\vc{q}}{\vt} \cdot \Grad \tvt \right] } \dt
	\label{w11}
\end{align}
for a.a. $0 < \tau < T$ whenever $\tvt \in C^1([0,T] \times \Ov{\Omega})$, $\tvt > 0$, $\tvt|_{\partial \Omega} = \vtB$.

\subsection{Constitutive relations}

The constitutive relation imposed on the equation of state and the transport coefficients are motivated by \cite{FeNo6A}. Specifically, we suppose 
\begin{equation} \label{w12}
	p(\vr, \vt) = \vt^{\frac{5}{2}} P \left( \frac{\vr}{\vt^{\frac{3}{2}}  } \right) + \frac{a}{3} \vt^4,\ 
	e(\vr, \vt) = \frac{3}{2} \frac{\vt^{\frac{5}{2}} }{\vr} P \left( \frac{\vr}{\vt^{\frac{3}{2}}  } \right) + \frac{a}{\vr} \vt^4, \ a > 0,
	\end{equation}
where $P \in C^1[0,\infty)$ satisfies
\begin{equation} \label{w13}
P(0) = 0,\ P'(Z) > 0 \ \mbox{for}\ Z \geq 0,\ 0 < \frac{ \frac{5}{3} P(Z) - P'(Z) Z }{Z} \leq c \ \mbox{for}\ Z > 0.
\end{equation} 	
This implies, in particular, that $Z \mapsto P(Z)/ Z^{\frac{5}{3}}$ is decreasing, and we suppose 
\begin{equation} \label{w14}
	\lim_{Z \to \infty} \frac{ P(Z) }{Z^{\frac{5}{3}}} = p_\infty > 0.
\end{equation}

Accordingly, the entropy $s$ takes the form 
\begin{equation} \label{w15}
s(\vr, \vt) = \mathcal{S} \left( \frac{\vr}{\vt^{\frac{3}{2}} } \right) + \frac{4a}{3} \frac{\vt^3}{\vr}, 
\end{equation}
where 
\begin{equation} \label{W15bis}
\mathcal{S}'(Z) = -\frac{3}{2} \frac{ \frac{5}{3} P(Z) - P'(Z) Z }{Z^2}.
\end{equation}
The reader may consult \cite[Chapter 1]{FeNo6A} for the physical background of the above constitutive theory. 

As for the transport coefficients, we suppose that $\mu$, $\eta$, and $\kappa$ are continuously differentiable functions of the temperature $\vt$ satisfying
\begin{align}
	0 < \underline{\mu} \left(1 + \vt^\Lambda \right) &\leq \mu(\vt) \leq \Ov{\mu} \left( 1 + \vt^\Lambda \right),\ 
	|\mu'(\vt)| \leq c \ \mbox{for all}\ \vt \geq 0,\ \frac{1}{2} \leq \Lambda \leq 1, \br
	0 &\leq  \eta(\vt) \leq \Ov{\eta} \left( 1 + \vt^\Lambda \right), \br 
	0 < \underline{\kappa} \left(1 + \vt^3 \right) &\leq  \kappa(\vt) \leq \Ov{\kappa} \left( 1 + \vt^3 \right).
	\label{w16}
	\end{align}
The property that both $\mu$ and $\kappa$ are unbounded for $\vt \to \infty$ is characteristic for gases, see Becker \cite{BE}.

\subsection{Weak solutions} 

We are ready to introduce the concept of \emph{weak solution} to the Navier--Stokes--Fourier system. 

\begin{Definition}[Weak solution] \label{wD1}
	We say that a trio $(\vr, \vt, \vu)$ is a \emph{weak solution} of the Navier--Stokes--Fourier system \eqref{i1}--\eqref{i6} in $(0,T) \times \Omega$ with the initial data 
	\begin{equation} \label{w17}
		\vr(0, \cdot) = \vr_0,\ (\vr \vu)(0, \cdot) = \vm_0,\ (\vr s(\vr, \vt))(0, \cdot) = S_0,
		\end{equation}
	if the following holds:
	\begin{itemize} 
		\item {\bf Conservation of mass.}
		\[
		\vr \in C_{\rm weak}([0,T]; L^{\frac{5}{3}}(\Omega)) \cap L^{\frac{5}{3}}((0,T) \times \partial \Omega; |\vuB \cdot \vc{n}|^+ (\dt \otimes \D \sigma_x) ),\ \vr \geq 0 \ 
		\mbox{a.a.};
		\]
		\[
		\vr \vu \in C_{\rm weak} ([0,T]; L^{\frac{5}{4}}(\Omega; R^d));
		\]
		the weak formulation of the equation of continuity \eqref{w1} is satisfied 
		with $\vr(0, \cdot) = \vr_0$
		for any $0 \leq \tau \leq T$, and any test function $\varphi \in C^1([0,T] \times \Ov{\Omega})$.
		\item {\bf Balance of momentum.} 
		\[
		\vu \in L^r(0,T; W^{1,r}(\Omega; R^d)),\ r = \frac{8}{5 - \Lambda}, \ (\vu - \vuB) \in L^r(0,T; W^{1,r}_0(\Omega; R^d));
		\]
		the weak formulation of the momentum balance \eqref{w3} holds with $(\vr \vu)(0, \cdot) = \vm_0$ for any $0 \leq \tau \leq T$, and any 
		$\bfphi \in C^1_c([0,T] \times \Omega; R^d)$.
		\item {\bf Entropy inequality.}
		\begin{align}
		\vt &\in L^\infty(0,T; L^4(\Omega)) \cap L^2(0,T; W^{1,2}(\Omega)),\ 
		(\vt - \vtB) \in L^2(0,T; W^{1,2}_0(\Omega)),\ \vt > 0 \ \mbox{a.a.}, \br 
		\log(\vt) &\in L^2((0,T); W^{1,2}(\Omega));
	\nonumber	
	\end{align}
the weak formulation of the entropy inequality \eqref{w5} is satisfied with $(\vr s)(0, \cdot) = S_0$ for any $\varphi \in C^1_c([0,T) \times \Omega)$, $\varphi \geq 0$.
\item {\bf Ballistic energy balance.}
The inequality \eqref{w10} holds for any 
\[
\tvt \in C^1([0,T] \times \Ov{\Omega}),\ \tvt > 0 ,\ \tvt|_{\partial \Omega} = \vtB,
\]
and any $\psi \in C^1_c[0,T)$, $\psi \geq 0$.

		\end{itemize}
	
	\end{Definition}

In the next section, we show that the weak solutions introduced in Definition \ref{wD1} comply with the weak--strong uniqueness principle and therefore represent a suitable generalization of classical solutions 

\section{Relative energy and weak--strong uniqueness principle}
\label{WS}

The proof of the weak--strong uniqueness principle is in the same spirit as in \cite{FeiNov20}. We introduce the \emph{relative energy} and use it as a Bregman distance between the strong and weak solution.

\subsection{Relative energy}

The relative energy is defined in the same way as in \cite{FeiNov20}:
\begin{align}
E &\left( \vr, \vt, \vu \Big| \tvr , \tvt, \tvu \right) \br &= \frac{1}{2}\vr |\vu - \tvu|^2 + \vr e - \tvt \Big(\vr s - \tvr s(\tvr, \tvt) \Big)- 
\Big( e(\tvr, \tvt) - \tvt s(\tvr, \tvt) + \frac{p(\tvr, \tvt)}{\tvr} \Big)
(\vr - \tvr) - \tvr e (\tvr, \tvt) \br &= 
 \frac{1}{2}\vr |\vu - \tvu|^2 + \vr e - \tvt \vr s - \Big( e(\tvr, \tvt) - \tvt s(\tvr, \tvt) + \frac{p(\tvr, \tvt)}{\tvr} \Big) \vr + p(\tvr, \tvt)
\label{ws1}
\end{align}
As observed in \cite{FeiNov20}, the relative energy interpreted in terms of the conservative entropy variables $(\vr, S = \vr s, \vm = \vr \vu)$ represents a Bregman distance associated to the energy functional 
\[
E(\vr, S, \vm) = \frac{1}{2} \frac{|\vm|^2}{\vr} + \vr e(\vr, S).
\]
Indeed it follows from our hypotheses concerning the form of the equation of state that $e$, $p$ satisfy the \emph{hypothesis of thermodynamic stability} 
\[
\frac{\partial p(\vr, \vt)}{\partial \vr} > 0,\ \frac{\partial e(\vr, \vt)}{\partial \vt} > 0,
\]
which in turn yields convexity of the internal energy $\vr e(\vr,S)$ with respect to the variables $(\vr,S)$. In addition,
\begin{equation} \label{ws2}
	\frac{\partial (\vr e(\vr, S))}{\partial \vr} = e - \vt s + \frac{p}{\vr},\ 
	\frac{\partial (\vr e(\vr, S))}{\partial S} = \vt.
	\end{equation} 
Thus the relative energy expressed in the conservative entropy variable may be interpreted as 
\[
E \left( \vr, S, \vm \Big| \tvr , \tS, \tvm \right) = E(\vr, S, \vm) - \left< \partial E(\tvr, \tS, \tvm) ; (\vr - \tvr, S - \tS, \vm - \tvm) \right> - E(\tvr, \tS, \tvm).
\]

\subsection{Relative energy inequality}

Our goal is to describe the time evolution of the relative energy 
\[
t \mapsto \intO{ E \left( \vr, \vt, \vu \Big| \tvr , \tvt, \tvu \right) (t, \cdot) },
\]
where $(\vr, \vt, \vu)$ is a weak solution of the Navier--Stokes--Fourier system and $(\tvr, \tvt, \tvu)$ is an arbitrary trio of ``test'' functions satisfying
\begin{align} 
	\tvr &\in C^1([0,T] \times \Ov{\Omega}),\ \inf \tvr > 0,\ \tvt \in C^1([0,T] \times \Ov{\Omega}),\ \inf \tvt > 0,\ 
	\tvt|_{\partial \Omega} = \vtB,\br \tvu &\in C^1([0,T] \times \Ov{\Omega}; R^d),\ \tvu|_{\partial \Omega} = \vuB.  
\label{ws3}	
\end{align}

Going back to \eqref{ws1} we observe that 
\begin{align}
\intO{ E \left( \vr, \vt, \vu \Big| \tvr , \tvt, \tvu \right) } &=  \underbrace{\intO{ \left[ \frac{1}{2}\vr |\vu - \tvu|^2 + \vr e - \tvt \vr s \right]}}_{\rm ballistic \ energy} \br &-  
\intO{ \Big( e(\tvr, \tvt) - \tvt s(\tvr, \tvt) + \frac{p(\tvr, \tvt)}{\tvr} \Big) \vr } +  \intO{p(\tvr, \tvt)},
\nonumber
\end{align}
where the time evolution of the ballistic energy is given by \eqref{w11}. Indeed the inequality \eqref{w11} is in fact independent of the specific extension of the boundary velocity $\vuB$ inside $\Omega$. 
As for the integral 
\[
\intO{ \Big( e(\tvr, \tvt) - \tvt s(\tvr, \tvt) + \frac{p(\tvr, \tvt)}{\tvr} \Big) \vr },
\]
it can be evaluated by using  $\Big( e(\tvr, \tvt) - \tvt s(\tvr, \tvt) + \frac{p(\tvr, \tvt)}{\tvr} \Big)$ as a test function in the weak form of the equation of continuity \eqref{w5}. 

Summarizing, we deduce the \emph{relative energy inequality} in the form
 \begin{align}
	&\left[ \intO{ E \left(\vr, \vt, \vu \Big| \tvr, \tvt, \tvu \right) } \right]_{t = 0}^{t = \tau} \br 
	&+ \int_0^\tau \intgi{ \left[ \vrB e(\vrB, \vtB) - \vtB \vrB s (\vrB, \vtB) - \left( e(\tvr, \vtB) - \vtB s(\tvr, \vtB) + \frac{p(\tvr, \vtB)}{\tvr} \right) \vrB  \right] } \dt
	\br 
	&+ \int_0^\tau \intgo{ \left[ \vr e(\vr, \vtB) - \vtB \vr s (\vr, \vtB) - \left( e(\tvr, \vtB) - \vtB s(\tvr, \vtB) + \frac{p(\tvr, \vtB)}{\tvr} \right) \vr  \right] } \dt \br 
	&+ \int_0^\tau \intO{ \frac{\tvt}{\vt} \left( \mathbb{S} (\vt, \Ds \vu) : \Ds \vu - \frac{\vc{q}(\vt, \Grad \vt)}{\vt} \right) } \dt \br 
	&\leq  \int_0^\tau \intO{ \frac{\vr}{\tvr} (\vu - \tvu ) \cdot \Grad p(\tvr, \tvt) } \dt \br 
	&- \int_0^\tau \intO{ \left( \vr (s - s(\tvr, \tvt)) \partial_t \tvt + \vr (s - s(\tvr, \tvt)) \vu \cdot \Grad \tvt + 
		\left( \frac{\vc{q}(\vt, \Grad \vt)}{\vt} \right) \cdot \Grad \tvt \right) } \dt \br 
	&- \int_0^\tau \intO{ \Big[ \vr (\vu - \tvu) \otimes (\vu - \tvu) + p(\vr, \vt) \mathbb{I} - \mathbb{S}(\vt, \Ds \vu) \Big] : \Ds \tvu } \dt \br 
	&+ \int_0^\tau \intO{ \vr \left[ \vc{g} - \partial_t \tvu - (\tvu \cdot \Grad) \tvu - \frac{1}{\tvr} \Grad p(\tvr, \tvt) \right] \cdot (\vu - \tvu) } \dt \br 
	&+ \int_0^\tau \intO{ \left[ \left( 1 - \frac{\vr}{\tvr} \right) \partial_t p(\tvr, \tvt) - \frac{\vr}{\tvr} \vu \cdot \Grad p(\tvr, \tvt) \right] } \dt
	\label{ws4}
	\end{align}
for a.a. $0 < \tau < T$, for any weak solution $(\vr, \vt, \vu)$ of the Navier--Stokes--Fourier system in the sense of Definition \ref{wD1} and any trio of test functions 
$(\tvr, \tvt, \tvu)$ belonging to the class \eqref{ws3}.

\subsection{Weak strong uniqueness}

The obvious idea how to compare a weak and strong solution is to use the strong solution $(\tvr, \tvt, \tvu)$ as the trio of test functions in the relative energy inequality \eqref{ws4}.
Keeping in mind that $(\vr, \vt, \vu)$ and $(\tvr, \tvt, \tvu)$ share the same initial and boundary data as well as the driving force $\vc{g}$, the inequality \eqref{ws4} simplifies to
\begin{align}
	&\intO{ E \left(\vr, \vt, \vu \Big| \tvr, \tvt, \tvu \right) (\tau, \cdot) }  - \int_0^\tau \intgi{ p(\vrB, \vtB)  } \dt
	\br 
	&+ \int_0^\tau \intgo{ \left[ \vr e(\vr, \vtB) - \vtB \vr s (\vr, \vtB) - \left( e(\tvr, \vtB) - \vtB s(\tvr, \vtB) + \frac{p(\tvr, \vtB)}{\tvr} \right) \vr  \right] } \dt \br 
	&+ \int_0^\tau \intO{ \frac{\tvt}{\vt} \left( \mathbb{S} (\vt, \Ds \vu) : \Ds \vu - \frac{\vc{q}(\vt, \Grad \vt) \cdot \Grad \vt}{\vt} \right) } \dt \br 
	&\leq  \int_0^\tau \intO{ \frac{\vr}{\tvr} (\vu - \tvu ) \cdot \Grad p(\tvr, \tvt) } \dt \br 
	&- \int_0^\tau \intO{ \left( \vr (s - s(\tvr, \tvt)) \partial_t \tvt + \vr (s - s(\tvr, \tvt)) \vu \cdot \Grad \tvt + 
		\left( \frac{\vc{q}(\vt, \Grad \vt)}{\vt} \right) \cdot \Grad \tvt \right) } \dt \br 
	&- \int_0^\tau \intO{ \Big[ \vr (\vu - \tvu) \otimes (\vu - \tvu) + p(\vr, \vt) \mathbb{I} - \mathbb{S}(\vt, \Ds \vu) \Big] : \Ds \tvu } \dt \br 
	&+ \int_0^\tau \intO{ \frac{\vr}{\tvr} \Div \mathbb{S} (\tvt, \Ds \tvu ) \cdot (\tvu - \vu) } \dt \br 
	&+ \int_0^\tau \intO{ \left[ \left( 1 - \frac{\vr}{\tvr} \right) \partial_t p(\tvr, \tvt) - \frac{\vr}{\tvr} \vu \cdot \Grad p(\tvr, \tvt) \right] } \dt.
	\label{ws5}
\end{align}

Moreover, regrouping certain integrals in \eqref{ws5} we obtain 
\begin{align}
	&\intO{ E \left(\vr, \vt, \vu \Big| \tvr, \tvt, \tvu \right) (\tau, \cdot) }  - \int_0^\tau \intgi{ p(\vrB, \vtB)  } \dt
	\br 
	&+ \int_0^\tau \intgo{ \left[ \vr e(\vr, \vtB) - \vtB \vr s (\vr, \vtB) - \left( e(\tvr, \vtB) - \vtB s(\tvr, \vtB) + \frac{p(\tvr, \vtB)}{\tvr} \right) \vr  \right] } \dt \br 
	&+ \int_0^\tau \intO{ \left( \frac{\tvt}{\vt} - 1 \right) \mathbb{S} (\vt, \Ds \vu) : \Ds \vu  } \dt + 
	\int_0^\tau \intO{ \left( 1 - \frac{\tvt}{\vt} \right) \left( \frac{\vc{q}(\vt, \Grad \vt)\cdot \Grad \vt}{\vt}\right) }  \dt \br &+ 
	\int_0^\tau \intO{ \frac{ \vc{q}(\vt, \Grad \vt) \cdot \Grad (\tvt - \vt)}{\vt} } \dt \br
	&+ \int_0^\tau \intO{ \Big( \mathbb{S} (\tvt, \Ds \tvu ) - \mathbb{S} (\vt, \Ds \vu) \Big): \Ds (\tvu - \vu) } \dt \br 
		&\leq  \int_0^\tau \intO{ (\vu - \tvu ) \cdot \Grad p(\tvr, \tvt) } \dt \br 
	&- \int_0^\tau \intO{ \left( \tvr (s - s(\tvr, \tvt)) \partial_t \tvt + \tvr (s - s(\tvr, \tvt)) \tvu \cdot \Grad \tvt 
	 \right) } \dt \br 
	&- \int_0^\tau \intO{  p(\vr, \vt)  \Div \tvu } \dt \br 
	&+ \int_0^\tau \intO{ \left[ \left( 1 - \frac{\vr}{\tvr} \right) \partial_t p(\tvr, \tvt) - \frac{\vr}{\tvr} \vu \cdot \Grad p(\tvr, \tvt) \right] } \dt + 
	\int_0^\tau  \intO{ R_1 } \dt
	\label{ws6}
\end{align}
with the quadratic error 
\begin{align}
R_1 =  &	\vr (\tvu - \vu) \otimes (\vu - \tvu) : \Ds \tvu 
+ \left( \frac{\vr}{\tvr} - 1 \right) \Div \mathbb{S} (\tvt, \Ds \tvu) (\tvu - \vu) \br 
&+ \left( \frac{\vr}{\tvr} - 1 \right) (\vu - \tvu) \cdot \Grad p(\tvr, \tvt) + \vr (s - s(\tvr, \tvt)) (\tvu - \vu) \cdot \Grad \tvt \br 
&- (\vr - \tvr)(s - s(\tvr, \tvt))(\partial_t \tvt + \tvu \cdot \Grad \tvt).
\label{ws7}
\end{align}

\subsubsection{Pressure}

Obviously,
\begin{align}
- &\intO{ \tvu \cdot \Grad p(\tvr, \tvt) } = - \intgi{ p(\vrB, \vtB) } - \intgo{ p(\tvr, \vtB) } \br &+ \intO{ p(\tvr, \tvt) \Div \tvu }; 
\nonumber
\end{align}
whence \eqref{ws7} can be written as 
\begin{align}
	&\intO{ E \left(\vr, \vt, \vu \Big| \tvr, \tvt, \tvu \right) (\tau, \cdot) }  
	\br 
	&+ \int_0^\tau \intgo{ \left[ \vr e(\vr, \vtB) - \vtB \vr s (\vr, \vtB) - \left( e(\tvr, \vtB) - \vtB s(\tvr, \vtB) + \frac{p(\tvr, \vtB)}{\tvr} \right) \vr  \right] } \dt \br
	&+ \int_0^\tau \intgo{ p(\tvr, \vtB) } \br 
	&+ \int_0^\tau \intO{ \left( \frac{\tvt}{\vt} - 1 \right) \mathbb{S} (\vt, \Ds \vu) : \Ds \vu  } \dt + 
	\int_0^\tau \intO{ \left( 1 - \frac{\tvt}{\vt} \right) \left( \frac{\vc{q}(\vt, \Grad \vt)\cdot \Grad \vt}{\vt}\right) }  \dt \br &+ 
	\int_0^\tau \intO{ \frac{ \vc{q}(\vt, \Grad \vt) \cdot \Grad (\tvt - \vt)}{\vt} } \dt \br
	&+ \int_0^\tau \intO{ \Big( \mathbb{S} (\tvt, \Ds \tvu ) - \mathbb{S} (\vt, \Ds \vu) \Big): \Ds (\tvu - \vu) } \dt \br 
	&\leq 
	- \int_0^\tau \intO{ \left( \tvr (s - s(\tvr, \tvt)) \partial_t \tvt + \tvr (s - s(\tvr, \tvt)) \tvu \cdot \Grad \tvt 
		\right) } \dt \br 
	&+ \int_0^\tau \intO{(p( \tvr, \tvt) -   p(\vr, \vt))  \Div \tvu } \dt \br 
	&+ \int_0^\tau \intO{ \left[ \left( 1 - \frac{\vr}{\tvr} \right) \Big( \partial_t p(\tvr, \tvt) -  \tvu \cdot \Grad p(\tvr, \tvt) \Big) \right] } \dt + 
	\int_0^\tau \intO{ R_2 } \dt,
	\label{ws8}
\end{align}
with
\begin{equation} \label{ws9}
R_2 = R_1 + \left(1 - \frac{\vr}{\tvr} \right) (\vu - \tvu) \cdot \Grad p(\tvr, \tvt) \,.
\end{equation}

Next, we use the identity
\[
\frac{\partial s (\tvr, \tvt) }{\partial \vr} = - \frac{1}{\tvr^2} \frac{\partial p(\tvr, \tvt) }{\partial \vt},
\]
to compute
\begin{align}
	\Big( 1 - &\frac{\vr}{\tvr} \Big) \left( \partial_t p(\tvr, \tvt) + \tvu \cdot
	\Grad p(\tvr, \tvt)  \right) + \Div \tvu \Big( p(\tvr, \tvt) - p(\vr, \vt) \Big)\br
	= &\Div \tvu \Big( p(\tvr, \tvt)  -\frac{\partial p(\tvr, \tvt) }{\partial \vr} (\tvr -\vr ) -\frac{\partial p(\tvr, \tvt)}{\partial \vt} (\tvt -\vt )- p(\vr, \vt) \Big)
	\br
	+ &\left( 1 - \frac{\vr}{\tvr} \right) \frac{\partial p (\tvr, \tvt) }{\partial \vt} \Big( \partial_t \tvt + \tvu \cdot \Grad \tvt \Big) -
	\frac{ \tvt - \vt }{ \tvr } \frac{\partial p (\tvr, \tvt) }{\partial \vt} \Big( \partial_t \tvr + \tvu \cdot \Grad \tvr \Big)
	\br
	= &\Div \tvu \Big( p(\tvr, \tvt)  -\frac{\partial p(\tvr, \tvt)}{\partial \vr} (\tvr -\vr ) -\frac{\partial p(\tvr, \tvt)}{\partial \vt} (\tvt -\vt )- p(\vr, \vt) \Big)
	\br
	- &\tvr (\tvr - \vr ) \frac{\partial s (\tvr, \tvt) }{\partial \vr} \Big( \partial_t \tvt + \tvu \cdot \Grad \tvt \Big)
	+  \tvr( \tvt - \vt) \frac{\partial s (\tvr, \tvt) }{\partial \vr} \Big( \partial_t \tvr + \tvu \cdot \Grad \tvr \Big)
	\br
	= &\Div \tvu \Big( p(\tvr, \tvt)  -\frac{\partial p(\tvr, \tvt)|}{\partial \vr} (\tvr -\vr ) -\frac{\partial p(\tvr, \tvt)}{\partial \vt} (\tvt -\vt )- p(\vr, \vt) \Big)
	\br
	- &\tvr (\tvr - \vr ) \frac{\partial s (\tvr, \tvt) }{\partial \vr} \Big( \partial_t \tvt + \tvu \cdot \Grad \tvt \Big)
	-  \tvr( \tvt - \vt) \frac{\partial s (\tvr, \tvt) }{\partial \vt} \Big( \partial_t \tvt + \tvu \cdot \Grad \tvt \Big)
	\br
	- &(\tvt - \vt) \Div \left( \frac{\widetilde{\vc{q}} } {\tvt} \right) +
	\left( 1 - \frac{\vt}{\tvt} \right) \left( \widetilde {\mathbb{S}} : \Grad \tvu -
	\frac{ \widetilde{\vc{q} } \cdot \Grad \tvt }{\tvt} \right).
	\nonumber
\end{align}

Consequently, \eqref{ws8} gives rise to 
\begin{align}
	&\intO{ E \left(\vr, \vt, \vu \Big| \tvr, \tvt, \tvu \right) (\tau, \cdot) }  
	\br 
	&+ \int_0^\tau \intgo{ \left[ \vr e(\vr, \vtB) - \vtB \vr s (\vr, \vtB) - \left( e(\tvr, \vtB) - \vtB s(\tvr, \vtB) + \frac{p(\tvr, \vtB)}{\tvr} \right) \vr  \right] } \dt \br
	&+ \int_0^\tau \intgo{ p(\tvr, \vtB) } \br 
	&+ \int_0^\tau \intO{ \left( \frac{\tvt}{\vt} - 1 \right) \mathbb{S} (\vt, \Ds \vu) : \Ds \vu  } \dt + 
	\int_0^\tau \intO{ \left( 1 - \frac{\tvt}{\vt} \right) \left( \frac{\vc{q}(\vt, \Grad \vt)\cdot \Grad \vt}{\vt}\right) }  \dt \br &+ 
	\int_0^\tau \intO{ \frac{ \vc{q}(\vt, \Grad \vt) \cdot \Grad (\tvt - \vt)}{\vt} } \dt \br
	&+ \int_0^\tau \intO{ \Big( \mathbb{S} (\tvt, \Ds \tvu ) - \mathbb{S} (\vt, \Ds \vu) \Big): \Ds (\tvu - \vu) } \dt \br 
	&\leq - \int_0^\tau \intO{ (\tvt - \vt) \Div \left( \frac{\vc{q}(\tvt, \Grad \tvt) }{\tvt} \right)} \dt \br &+ \int_0^\tau \intO{ 
			\left(1 - \frac{\vt}{\tvt} \right) \left( \mathbb{S} (\tvt, \Ds \tvu) : \Ds \tvu - \frac{\vc{q} (\tvt, \Grad \tvt) }{\tvt} \right) } \dt + 
	\int_0^\tau \intO{ R_3 } \dt
	\label{ws10}
\end{align}
with the quadratic error
\begin{align}
R_3 &= R_2 + \Div \tvu \Big( p(\tvr, \tvt)  -\frac{\partial p(\tvr, \tvt)}{\partial \vr} (\tvr -\vr ) -\frac{\partial p(\tvr, \tvt)}{\partial \vt} (\tvt -\vt )- p(\vr, \vt) \Big) \br 
&+ \tvr \left( s(\tvr, \tvt) - \frac{\partial s(\tvr, \tvt)}{\partial \vr}(\tvr - \vr) - \frac{\partial s(\tvr, \tvt)}{\partial \vt}(\tvt - \vt) - s \right) (\partial_t \tvt + \tvu \cdot \Grad \tvt ).
\label{ws11}
\end{align}

Finally, using the fact that $\vt - \tvt$ vanishes on the boundary, we may integrate by parts
\[
- \intO{ (\tvt - \vt) \Div \left( \frac{\vc{q}(\tvt, \Grad \tvt) }{\tvt} \right)} = 
\intO{ \left( \frac{\vc{q}(\tvt, \Grad \tvt) }{\tvt} \right) \cdot \Grad (\tvt - \vt) }
\]
to conclude 	
\begin{align}
	&\intO{ E \left(\vr, \vt, \vu \Big| \tvr, \tvt, \tvu \right) (\tau, \cdot) }  
	\br 
	&+ \int_0^\tau \intgo{ \left[ \vr e(\vr, \vtB) - \vtB \vr s (\vr, \vtB) - \left( e(\tvr, \vtB) - \vtB s(\tvr, \vtB) + \frac{p(\tvr, \vtB)}{\tvr} \right) \vr  \right] } \dt \br
	&+ \int_0^\tau \intgo{ p(\tvr, \vtB) } \br 
	&+ \int_0^\tau \intO{ \left( \frac{\tvt}{\vt} - 1 \right) \mathbb{S} (\vt, \Ds \vu) : \Ds \vu  } \dt 
	+ \int_0^\tau \intO{ 
		\left(\frac{\vt}{\tvt} - 1 \right)  \mathbb{S} (\tvt, \Ds \tvu) : \Ds \tvu } \dt
	\br
	&+\int_0^\tau \intO{ \left( 1 - \frac{\tvt}{\vt} \right) \left( \frac{\vc{q}(\vt, \Grad \vt)\cdot \Grad \vt}{\vt}\right) }  \dt
	+ \int_0^\tau \intO{ 
		\left(1 - \frac{\vt}{\tvt} \right) \frac{\vc{q} (\tvt, \Grad \tvt) }{\tvt} } \dt
	\br &+ 
	\int_0^\tau \intO{ \left( \frac{ \vc{q}(\vt, \Grad \vt)}{\vt} - \frac{ \vc{q}(\tvt, \Grad \tvt)}{\tvt} \right)  \cdot \Grad (\tvt - \vt)}  \dt \br
	&+ \int_0^\tau \intO{ \Big( \mathbb{S} (\tvt, \Ds \tvu ) - \mathbb{S} (\vt, \Ds \vu) \Big): \Ds (\tvu - \vu) } \dt \leq  
	\int_0^\tau \intO{ R_3 } \dt
	\label{ws12}
\end{align}	
with the quadratic error $R_3$ determined successively via \eqref{ws7}, \eqref{ws9}, and \eqref{ws11}. 

The final inequality \eqref{ws12} as well as the quadratic error \eqref{ws11} are exactly the same (modulo the boundary flux not present in \eqref{ws12}) as in 
\cite[Section 4, formula (4.5)]{FeiNov20}. Thus repeating the arguments of \cite{FeiNov20} we deduce the desired conclusion 
\[
\intO{ E \left(\vr, \vt, \vu \Big| \tvr, \tvt, \tvu \right) (\tau, \cdot) } = 0 \ \mbox{for a.a.}\ \tau \in (0,T).
\] 
 
We have proved the following result. 

\begin{Theorem}[Weak--strong uniqueness principle]  \label{wsT1}
	Let the pressure $p = p(\vr, \vt)$, the internal energy $e = e(\vr, \vt)$, and the entropy $s = s(\vr, \vt)$ satisfy the hypotheses \eqref{w12}--\eqref{w15}. Let the transport coefficients 
	$\mu = \mu(\vt)$, $\eta = \eta(\vt)$, $\kappa = \kappa(\vt)$ be continuously differentiable functions of $\vt$ satisfying \eqref{w16}. Let $(\vr, \vt, \vu)$ be a weak solution 
	of the Navier--Stokes--Fourier system in $(0,T) \times \Omega$ in the sense of Definition \ref{wD1}. Suppose that 
	the same problem with the same initial data $\vr_0$, $\vm_0$, $S_0$, the boundary data $\vrB$, $\vtB$, $\vuB$, and the driving force $\vc{g}$ admits a strong solution 
	$(\tvr, \tvt, \tvu)$ in the class 
	\[
	\tvr, \ \tvt \in C^1([0,T] \times \Ov{\Omega}),\ 
	\tvu \in C^1([0,T] \times \Ov{\Omega}),\ \Grad^2 \tvt \in BC((0,T) \times \Omega; R^d),\ 
	\Grad^2 \tvu \in BC((0,T) \times \Omega; R^{d \times d}).
	 \]
	 
	 Then 
	 \[
	 \vr = \tvr,\ \vt = \tvt,\ \vu = \tvu \ \mbox{in}\ (0,T) \times \Omega.
	 \]
	
	\end{Theorem}

\section{Existence of weak solutions}
\label{E}

Our ultimate goal is to establish the \emph{existence} of a weak solution for given data on an arbitrary time interval $(0,T)$. 
We address in detail two principal issues:
\begin{itemize}
	\item {\it a priori} bounds that guarantee boundedness of the ballistic energy as well as the dissipation term in the inequality \eqref{w9};
	\item a suitable approximation scheme to construct the weak solutions.
	\end{itemize}
Having established this, we omit the proof of convergence of approximate solutions as 
it is basically the same as in \cite{FeiNov20}. As a matter of fact, the present setting is even easier than in \cite{FeiNov20} as the temperature is given on the whole domain boundary. 
To avoid technicalities, we assume that $\Omega \subset R^d$ is a smooth bounded domain of class at least $C^3$, and that the boundary data $\vrB$, $\vuB$, $\vtB$ are at least twice continuously differentiable.
In particular, the outer normal is well defined. 

\subsection{{\it A priori} bounds}

In comparison with \cite{FeNo6A} or \cite{FeiNov20}, the {\it a priori} bounds for the Dirichlet problem are more difficult to obtain as the heat flux through the boundary is not controlled. 
The relevant estimates are derived from the ballistic energy inequality \eqref{w9} evaluated for a \emph{suitable} temperature $\tvt$. To simplify, we shall assume that the viscosity coefficients 
satisfy hypothesis \eqref{w16} with $\Lambda = 1$. Consequently, by virtue of Korn--Poincar\' e inequality, 
\begin{equation} \label{AP1}
	\intO{ \frac{\tvt}{\vt} \mathbb{S}(\vt; \Ds \vu) : \Ds \vu } \ageq \inf_{\Omega}\{ \tvt \} \left( \| \vu \|^2_{W^{1,2}(\Omega)} - c(\vuB) \right).
	\end{equation}
By the same token, we get 
\[
\left| \intO{ \mathbb{S} (\vr ,\Ds \vu) : \Ds \vuB } \right| \leq \omega \| \Ds \vu \|^2_{L^2(\Omega; R^{d \times d})} + c(\vuB) (1 + \intO{ \vt^2 }). 
\]
Thus we may rewrite the ballistic energy inequality \eqref{w11} in the form 
\begin{align}  
	&\intO{ \left( \frac{1}{2} \vr |\vu - \vuB|^2 + \vr e - \tvt \vr s \right) (\tau, \cdot) }  \br  &+ 
	\int_0^\tau \intgo{ \Big[ \vr e(\vr, \vtB) - \vtB \vrB s(\vr, \vtB) \Big] } \dt \br &+ 
	\int_0^\tau \inf_{\Omega} \{ \tvt \} \intO{ \left( \| \vu \|^2_{W^{1,2}(\Omega; R^d)} + \frac{\kappa (\vt) |\Grad \vt|^2 }{\vt^2} \right) } \dt \br
	&\leq 
	\intO{ \left( \frac{1}{2} \vr_0 |\vu_0 - \vuB|^2 + \vr_0 e(\vr_0, \vt_0)  - \tvt \vr_0 s (\vr_0, \vt_0) \right)  } \br 
	& + c(\vrB, \vtB, \vuB, \vc{g} ) \Big[ 1 +
	\int_0^\tau  \intO{ \left( \frac{1}{2} \vr |\vu - \vuB|^2 + \vr e - \tvt \vr s \right)  } \dt 
	\br
	&- \int_0^\tau \intO{ \left[ \vr s \left( \partial_t \tvt + \vu \cdot \Grad \tvt \right) + \frac{\vc{q}}{\vt} \cdot \Grad \tvt \right] } \dt \Big]\,.
	\label{AP2}
\end{align}

It remains to control the last integral in \eqref{AP2}. To this end, we first fix the extension $\tvt$ to be the unique solution of the Laplace equation 
\begin{equation} \label{AP3}
	\Del \tvt (\tau, \cdot) = 0 \ \mbox{in}\ \Omega,\ \tvt(\tau, \cdot)|_{\partial \Omega} = \vtB \ \mbox{for any}\ \tau \in [0,T].
	\end{equation}
It follows from the standard maximum principle for harmonic functions that 
\[
\min_{[0,T] \times \partial \Omega} \vtB \leq \tvt (t,x) \leq \max_{[0,T] \times \partial \Omega} \vtB \ \mbox{for any}\ (t,x) \in (0,T) \times \Omega.
\]
Let us denote this particular extension as $\vtB$. 

A simple integration by parts yields 
\[
- \intO{ \frac{\vc{q}}{\vt} \cdot \Grad \vtB } = \intO{ \frac{\kappa(\vt)}{\vt} \Grad \vt \cdot \Grad \vtB } = 
\intO{ \Grad K(\vt) \cdot \Grad \vtB } = \int_{\partial \Omega} K(\vtB) \Grad \vtB \cdot \vc{n},  
\] 
where $K'(\vt) = \frac{\kappa (\vt)}{\vt}$. 

Next, as $\vtB$ is continuously differentiable in time, we get
\[
- \intO{ \vr s \partial_t \vtB } \leq c_1 + c_2 \intO{  \left( \frac{1}{2} \vr |\vu - \vuB|^2 + \vr e - \vtB \vr s \right) } 
\ \mbox{for some}\ c_1, c_2 > 0.
\]
Summarizing we deduce from \eqref{AP2}: 
\begin{align}  
	&\intO{ \left( \frac{1}{2} \vr |\vu - \vuB|^2 + \vr e - \vtB \vr s \right) (\tau, \cdot) }  \br  &+ 
	\int_0^\tau \intgo{ \Big[ \vr e(\vr, \vtB) - \vtB \vrB s(\vr, \vtB) \Big] } \dt \br &+ 
	\inf_{[0,T] \times \partial \Omega}\{ \vtB \} \int_0^\tau  \intO{ \left( \| \vu \|^2_{W^{1,2}(\Omega; R^d)} + \frac{\kappa (\vt) |\Grad \vt|^2 }{\vt^2} \right) } \dt \br
	&\leq 
	\intO{ \left( \frac{1}{2} \vr_0 |\vu_0 - \vuB|^2 + \vr_0 e(\vr_0, \vt_0)  - \vtB \vr_0 s (\vr_0, \vt_0) \right)  } + c_1 \br 
	& + c_2
	\int_0^\tau  \intO{ \left( \frac{1}{2} \vr |\vu - \vuB|^2 + \vr e - \vtB \vr s \right)  } \dt 
	\br
	&- \int_0^\tau \intO{ \vr s \vu \cdot \Grad \vtB  } \dt  \,.
	\label{AP4}
\end{align}

As the entropy is given by \eqref{w15}, we have 
\begin{equation} \label{AP5}
- \intO{ \vr s \vu \cdot \Grad \vtB  } = - \intO{ \vr \mathcal{S} \left( \frac{\vr}{\vt^{\frac{3}{2} } } \right) \vu \cdot \Grad \vtB } - \frac{4a}{3} \intO{ \vt^3 \vu \cdot \Grad \vtB },
\end{equation}
where $\mathcal{S}$ is determined, modulo an additive constant, by \eqref{W15bis}. In addition to the hypotheses \eqref{w12}--\eqref{w14}, we suppose that $s$ satisfies 
the Third law of thermodynamics, specifically, 
\begin{equation} \label{AP6bis}
s(\vr, \vt) \to 0 \ \mbox{as}\ \vt \to 0 \ \mbox{for any}\ \vr > 0.
\end{equation}
As shown in \cite[Section 4, formula (4.6)]{FeiNov10}, this implies 
\[ 
\varrho \left\vert \mathcal{S} \left( \frac{\vr}{\vt^{\frac{3}{2} } } \right) \right \vert  \aleq \left( \vr + \vr |\log(\vr)| + \vr [ \log(\vt)]^+ \right) \ \mbox{for any} \ \vr \geq 0,\ \vt \geq 0.\]
Consequently, going back to \eqref{AP5} we obtain 
\begin{align}
&\left| \intO{ \vr \mathcal{S} \left( \frac{\vr}{\vt^{\frac{3}{2} } } \right) \vu \cdot \Grad \vtB } \right| \br &\aleq 
\left( \intO{ \vr |\vu|^2 } + \intO{\vr \mathcal{S}^2\left( \frac{\vr}{\vt^{\frac{3}{2} } } \right) } \right) \aleq 
\left( 1 + 	 \intO{ \vr |\vu|^2 } + \intO{ \vr^{\frac{5}{3} } } + \intO{\vt^4} \right)  \br 
&\leq c_3 +  c_4 \intO{ \left( \frac{1}{2} \vr |\vu - \vuB|^2 + \vr e - \vtB \vr s \right) },\ c_3, c_4 > 0.
	\label{AP6}
	\end{align}

In order to control the radiative component of the entropy flux, we have to strengthen our hypothesis concerning the growth of the heat conductivity coefficient $\kappa(\vt)$ for large values of $\vt$. 
In accordance with our previous agreement concerning the linear growth of the viscosity coefficients, we replace \eqref{w16} by
\begin{align}
	0 < \underline{\mu} \left(1 + \vt \right) &\leq \mu(\vt) \leq \Ov{\mu} \left( 1 + \vt \right),\ 
	|\mu'(\vt)| \leq c \ \mbox{for all}\ \vt \geq 0,\br
	0 &\leq  \eta(\vt) \leq \Ov{\eta} \left( 1 + \vt \right), \br 
	0 < \underline{\kappa} \left(1 + \vt^\beta \right) &\leq  \kappa(\vt) \leq \Ov{\kappa} \left( 1 + \vt^\beta \right),\ \beta > 6.
	\label{w16bis}
\end{align} 

Consequently, on the one hand, 
\begin{equation} \label{AP8}
\intO{ \frac{\kappa (\vt)|\Grad \vt|^2 }{\vt^2} } \ageq \intO{ \left( \frac{1}{\vt^2} + \vt^{\beta - 2} \right) |\Grad \vt|^2  },
\end{equation}
while, on the other hand, 
\begin{equation} \label{AP9}
\left| \intO{ \vt^3 \vu \cdot \Grad \vtB } \right| \leq \ep \| \vu \|^2_{W^{1,2}(\Omega; R^d)} + c(\ep) \| \vt^3 \|_{L^2(\Omega)}^2
\end{equation}
for any $\ep > 0$. Furthermore, 
\[
\| \vt^3 \|_{L^2(\Omega)}^2 = \int_{\vt \leq K} \vt^6 \ \dx +   \int_{\vt > K} \vt^6 \ \dx \leq 
|\Omega| K^6 + K^{6 - \beta} \intO{ \vt^\beta }.
\]
Next, by H\" older and Poincar\' e inequalities, 
\begin{align}
\intO{ \vt^\beta } &\aleq \| \vt^\beta \|_{L^3(\Omega)} = \| \vt^{\frac{\beta}{2}} \|^2_{L^6(\Omega)} \aleq 
\| \vt^{\frac{\beta}{2}} \|^2_{W^{1,2}(\Omega)} \aleq \| \Grad \vt^{\frac{\beta}{2}} \|^2_{L^2(\Omega)} + c(\vtB) 
\br &\aleq \intO{ \vt^{\beta - 2} |\Grad \vt|^2 } + c(\vtB).
\nonumber
\end{align}
As $\beta > 6$, we can first fix $\ep > 0$ small enough and then $K = K(\ep)$ large enough so that 
the integral \eqref{AP9} is dominate by the left--hand side of \eqref{AP4}. 

Summing up the previous estimates we may convert \eqref{AP4} in the desired relation 
\begin{align}  
	&\intO{ \left( \frac{1}{2} \vr |\vu - \vuB|^2 + \vr e - \tvt \vr s \right) (\tau, \cdot) }  \br  &+ 
	\int_0^\tau \intgo{ \Big[ \vr e(\vr, \vtB) - \vtB \vrB s(\vr, \vtB) \Big] } \dt \br &+ 
	\inf_{[0,T] \times \partial \Omega}\{ \vtB \} \int_0^\tau  \intO{ \left( \| \vu \|^2_{W^{1,2}(\Omega; R^d)} + \frac{\kappa (\vt) |\Grad \vt|^2 }{\vt^2} \right) } \dt \br
	&\leq 
	\intO{ \left( \frac{1}{2} \vr_0 |\vu_0 - \vuB|^2 + \vr_0 e(\vr_0, \vt_0)  - \tvt \vr_0 s (\vr_0, \vt_0) \right)  } + c_5 \br 
	& + c_6
	\int_0^\tau  \intO{ \left( \frac{1}{2} \vr |\vu - \vuB|^2 + \vr e - \tvt \vr s \right)  } \dt . 
	\label{AP10}
\end{align}
The required {\it a priori} bounds follow form \eqref{AP10} by means of Gronwall lemma. 

\begin{Remark}[Constant boundary temperature] \label{RAP1}
	Of course, the extra hypotheses \eqref{AP6bis}, \eqref{w16bis} are not necessary provided the boundary temperature $\vtB$ is a \emph{positive} constant
	or a function that depends only on $t$, in which case $\Grad \vtB = 0$. 
	
	\end{Remark}

\subsection{Approximation scheme}

First we choose a system of functions $\{ \vc{w}_j \}_{j=1}^\infty \subset \DC(\Omega; R^d)$ that form an orthonormal basis of the space $L^2(\Omega; R^d)$ and fix two positive parameters $\ep > 0$, $\delta > 0$.
We also introduce an auxiliary function
\[
[z]^-_n \in C^\infty(R),\ [z]^-_n = \left\{ \begin{array}{l} z \ \mbox{if}\ z < - \frac{1}{n}, \\ 
	\mbox{non--decreasing} \leq \min\{ z, 0 \} \ \mbox{if}\ - \frac{1}{n} \leq z \leq \frac{1}{n}, \\
	0 \ \mbox{if}\ z > \frac{1}{n}. \end{array} \right.
\] 
In other words, $[z]^-_n$ is a smooth approximation of the negative part $[z]^- = \min \{ 0, z \}$.

The basic level approximate solution is a trio of functions $(\vr, \vt, \vu)$, where 
\[
\vu = \vc{v} + \vuB,\ \vc{v} \in C^1([0,T]; X_n),\ X_n = {\rm span} \left\{ \vc{w}_j \ \Big| \ 1 \leq j \leq n \right\}.
\] 
The functions $(\vr, \vt, \vu)$ solve the following system of equations: 

\begin{itemize}
	\item {\bf Vanishing viscosity approximation of the equation of continuity.}
	\begin{align} 
	\partial_t \vr + \Div (\vr \vu) &= \ep \Del \vr \ \mbox{in}\ (0,T) \times \Omega, \br
	\ep \Grad \vr \cdot \vc{n} &= (\vr - \vrB) [\vuB \cdot \vc{n}]^-_n \ \mbox{in}\ (0,T) \times \partial \Omega, \br 
	\vr(0, \cdot) &= \vr_{0,\delta}.
		\label{E1}	
		\end{align}
	
	\item {\bf Galerkin approximation of the momentum equation.}	
	
\begin{align}
	\int_0^\tau &\intO{ \Big[ \vr \vu \cdot \partial_t \bfphi + \vr \vu \otimes \vu : \Grad \bfphi + \delta \left( \vr^\Gamma + \vr^2 \right) \Div \bfphi + p(\vr, \vt) \Div \bfphi \Big] } \dt \br &= 
	\int_0^\tau \intO{ \Big[ \mathbb{S}_\delta (\vt, \Ds \vu) : \Ds \bfphi - \vr \vc{g} \cdot \bfphi \Big] } \dt + \left[ \intO{ \vr \vu \cdot \bfphi } \right]_{t=0}^{t = \tau}, \ 
	\vr \vu (0, \cdot) = \vm_0,  
	\label{E2}
\end{align}
for any $0 \leq \tau \leq T$, $\bfphi \in C^1([0,T]; X_n)$. Here $\mathbb{S}_\delta$ denotes the viscous stress with the shear viscosity $\mu_\delta = \mu + \delta \vt$.

\item {\bf Approximate internal energy balance.}

\begin{align}
	\partial_t \Big[ \vr (e(\vr, \vt) + \delta \vt) \Big] &+ \Div \Big[ \vr (e(\vr, \vt) + \delta  \vt) \vu \Big] - \Div \left[\left(\delta \left( \vt^\Gamma + \frac{1}{\vt} \right)  +    \kappa (\vt) \right) \Grad \vt \right] \br & = \mathbb{S}_\delta : \Ds \vu - 
	p(\vr, \vt) \Div \vu \br &+ \ep \delta \left( \Gamma \vr^{\Gamma - 2} + 2 \right) |\Grad \vr|^2 + \frac{\delta}{\vt^2} - \ep \vt^5 \ \mbox{in}\ (0,T) \times \Omega, \br
	\vt&= \vtB \ \mbox{in}\ (0,T) \times \partial \Omega, \br
	\vt(0, \cdot) &= \vt_0.
	\label{E3}
	\end{align}	

\end{itemize}

Thus the scheme at the basic level is almost identical with that used in \cite{FeiNov20} except the boundary conditions in \eqref{E3}. Consequently, we can pass to the entropy as well as the ballistic energy balance as long as the temperature is strictly positive. Indeed it is legal to use $(\vu - \vuB)$ as a test function in the approximate momentum balance and to integrate the approximate internal energy equation to 
obtain an analogue of the total energy balance \eqref{w7} as well as the ballistic energy inequality \eqref{w11}. With the {\it a priori estimates} established in the preceding section at hand, the convergence of the 
approximate scheme can be shown similarly to \cite{FeiNov20}. We conclude by stating the global existence result. 

\begin{Theorem}[Global existence]\label{APT1}
	Let $\Omega \subset R^d$, $d = 2,3$ be a bounded domain of class at least $C^3$. Suppose that the boundary data $\vuB = \vuB(t,x)$, $\vrB = \vrB(t,x)$, and $\vtB = \vtB(t,x)$ are twice continuously differentiable, and
	\[
	\inf_{[0,T] \times \partial \Omega } \vrB > 0,\ \inf_{[0,T] \times \partial \Omega} \vtB > 0.
	\]
	Let the pressure $p = p(\vr, \vt)$, the internal energy $e = e(\vr, \vt)$, the entropy $s = s(\vr, \vt)$ satisfy the hypotheses \eqref{w12}--\eqref{W15bis}, and let the transport coefficients 
	$\mu = \mu(\vt)$, $\eta = \eta (\vt)$, and $\kappa = \kappa(\vt)$ satisfy \eqref{w16}. In addition, if $\vtB \ne {\rm const}$, we suppose \eqref{AP6bis} and \eqref{w16bis}. 
	
	Then for any $T > 0$, any initial data 
	\[
	\vr_0,\ \vt_0, \vu_0,\ \vr_0 > 0, \ \vt_0 > 0,\ \intO{ \left( \frac{1}{2} \vr_0 |\vu_0|^2 + \vr_0 e(\vr_0, \vt_0)  - \vtB \vr_0 s (\vr_0, \vt_0) \right)  } < \infty, 
	\]
	and any $\vc{g} \in L^\infty((0,T) \times \Omega; R^d)$, 
	there exists a weak solution $(\vr, \vt, \vu)$ of the Navier--Stokes--Fourier system in $(0,T) \times \Omega$ in the sense of Definition \ref{wD1}.
	
	\end{Theorem}

\def\cprime{$'$} \def\ocirc#1{\ifmmode\setbox0=\hbox{$#1$}\dimen0=\ht0
	\advance\dimen0 by1pt\rlap{\hbox to\wd0{\hss\raise\dimen0
			\hbox{\hskip.2em$\scriptscriptstyle\circ$}\hss}}#1\else {\accent"17 #1}\fi}

%\bibliography{citace}
%\bibliographystyle{plain}

\end{document}